\documentclass[10pt,draft]{amsart}
\usepackage{enumerate}

\newtheorem{theorem}{Theorem}[section]
\newtheorem{proposition}[theorem]{Proposition}
\newtheorem{corollary}[theorem]{Corollary}
\newtheorem{lemma}[theorem]{Lemma}
\theoremstyle{definition}
\newtheorem{remark}[theorem]{Remark}

\newtheorem{notation}[theorem]{Notation}

\newtheorem{condition}[theorem]{Condition}

\def\Diff{\Delta}

\def\QQ{\mathbb{Q}}
\def\NN{\mathbb{N}}
\def\N{\mathbb{N}}
\def\ZZ{\mathbb{Z}}
\def\Z{\mathbb{Z}}
\def\RR{\mathbb{R}}
\def\R{\mathbb{R}}

\def\IA{\mathcal{A}}

\def\sm{\setminus}

\def\xelm#1{\mathbf{\widetilde{#1}}}

\def\Id{\mathrm{Id}}
\def\rest{\hbox{\kern0.1em\rule[-0.4em]{0.4pt}{1.1em}\kern0.1em}}

\begin{document}

\title[Invariant decomposition of functions]{Invariant decomposition
of functions with respect to commuting invertible transformations}

\author[B. Farkas]{B\'alint Farkas}
\address[B. Farkas]{Technische Universit\"{a}t Darmstadt, Fachbereich Mathematik,
AG4\newline Schlo\ss{}gartenstra\ss{}e 7, D-64289, Darmstadt,
Germany} \email{farkas@mathematik.tu-darmstadt.de}

\author[V. Harangi]{Viktor Harangi}

\author[T. Keleti]{Tam\'as Keleti}

\address[V. Harangi, T. Keleti]{
 Department of Analysis, E\"otv\"os Lor\'and University
\newline P\'az\-m\'any P\'e\-ter s\'et\'any 1/c, H-1117 Budapest,
Hungary}
 \email{bizkit@cs.elte.hu}
 \email{elek@cs.elte.hu}

\author[Sz. Gy. R\'ev\'esz]{Szil\'ard Gy\"orgy R\'{e}v\'{e}sz}
\address[Sz. Gy. R\'ev\'esz]{
A. R\'enyi Institute of Mathematics, Hungarian Academy of Sciences
\newline Budapest, P.O.B. 127, 1364 Hungary.}
\email{revesz@renyi.hu}

\address{and}
\address{Institut Henri Poincar\'e
\newline 11 rue Pierre et Marie Curie, 75005 Paris, France} \email{revesz@ihp.jussieu.fr}

\begin{abstract}
Consider $a_1,\dots,a_n\in\RR$ arbitrary elements. We characterize
those functions $f:\RR\to\RR$  that decompose into the sum of
$a_j$-periodic functions, i.e., $f=f_1+\cdots+f_n$ with
$\Diff_{a_j}f(x):=f(x+a_j)-f(x)=0$. We show that $f$ has such a
decomposition if and only if for all partitions $B_1\cup
B_2\cup\cdots \cup B_N=\{a_1,\dots,a_n\}$  with $B_j$ consisting of
commensurable elements with least common multiples $b_j$ one has
$\Diff_{{b_1}}\dots \Diff_{{b_N}}f=0$.

Actually, we prove a more general result for periodic decompositions
of functions $f:\mathcal{A}\to \RR$ defined on an Abelian group
$\mathcal{A}$, in fact, we even consider invariant decompositions of
functions $f:A\to \RR$ with respect to commuting, invertible
self-mappings of some abstract set $A$.

We also extend our results to functions between torsion free Abelian
groups. As a corollary we also obtain that on a torsion free Abelian
group the existence of a real valued periodic decomposition of an
integer valued function implies the existence of an integer valued
periodic decomposition with the same periods.

\end{abstract}


\thanks{Supported in the framework of the Hungarian-Spanish
Scientific and Technological Governmental Cooperation, Project \#
E-38/04 and in the framework of the Hungarian-French Scientific and
Technological Governmental Cooperation, Project \# F-10/04.}

\thanks{The third author was supported by Hungarian Scientific
Foundation grants no.~F 43620 and T 49786.}

\thanks{This work was accomplished during the fourth author's
stay in Paris under his Marie Curie fellowship, contract \#
MEIF-CT-2005-022927.}

\maketitle


\section{Introduction}
 The starting point of this note is the following
observation. If we have $a_j$-periodic functions $f_j:\RR\to\RR$,
$j=1,\dots,n$, then the sum $f:=f_1+f_2+\cdots+f_n$ satisfies the
difference equation
\begin{equation}\label{eq:diffeq}
\Diff_{a_1}\Diff_{a_2}\cdots \Diff_{a_n}f=0,\quad\mbox{where $\quad
\Diff_{a_j}f(x)=f(x+a_j)-f(x)$.}
\end{equation}
The converse implication, i.e., that the above difference equation
would imply existence of a periodic decomposition, however fails
already in the simplest situation. For instance, take
$a_1=a_2=a\in\RR$ and $f=\Id:\RR\to\RR$, the identity. Then
$\Diff_{a_1}\Diff_{a_2}f=0$ holds but of course $f$ is not a sum of
 two $a$-periodic functions, as it is not $a$-periodic.

\par  There are two natural ways to overcome this. One might restrict
the whole question by requiring that both $f$ and the $f_j$s belong
to some function class $\mathcal{F}$, which is then said to have the
\emph{decomposition property}, if the existence of a periodic
decomposition in $\mathcal{F}$ is equivalent to the above difference
equation \eqref{eq:diffeq}. For example, it is known that the class
$B(\RR)$ of bounded functions \cite{laczkovich/revesz:1990}, the
class $BC(\RR)$ of bounded continuous functions
\cite{laczkovich/revesz:1989} or the class $UCB(\R)$ of uniformly
continuous bounded functions \cite{laczkovich/revesz:1989},
\cite{gajda:1992}, or more generally $UCB(\mathcal{A})$ for any
locally compact topological group $\mathcal{A}$ \cite{gajda:1992}
\emph{do have}, while the above example shows that $C(\RR)$ and
$\RR^\RR$ \emph{do not have} the decomposition property. Some
natural questions are still open, for example it is not known, to
the best of our knowledge, whether $BC(\mathcal{A})$ has the
decomposition property for any  locally compact topological group
$\mathcal{A}$.
 For more about the
decomposition property of function classes  see Kadets, Shumyatskiy
\cite{kadets/shumyatskiy:2000} \cite{kadets/shumyatskiy:2001},
Keleti \cite{keleti:1996,keleti:1997,keleti:1998} and Laczkovich,
R\'{e}v\'{e}sz \cite{laczkovich/revesz:1989,
laczkovich/revesz:1990}.

\vskip1ex The other possibility, which is actually our goal now,
instead of restricting to a particular function class, is to
complement the above difference equation with other conditions of
similar type, which then together will be sufficient and necessary
for the existence of periodic decompositions with given periods.
Suppose that $f$ has a periodic decomposition with periods
$a_1,\dots,a_n$ and let $B_1\cup\dots\cup B_N=\{a_1,\dots,a_n\}$ be
a partition such that in each $B_j$ the elements are commensurable
(that is, they have a common  (nonzero integer) multiple)
 with least common multiple $b_j$. Then, by summing up
for each $j$ the terms corresponding to the elements in $B_j$ we get
a periodic decomposition of $f$ with periods $b_1,\dots,b_N$. Thus
we must have $\Diff_{b_1}\Diff_{b_2}\cdots \Diff_{b_N}f=0$.
Therefore we see that if $f$ has a periodic decomposition with
periods $a_1,\dots,a_n$, then for any partition $B_1\cup\dots\cup
B_N=\{a_1,\dots,a_n\}$ such that in each $B_j$ the elements are
commensurable with least common multiple $b_j$, we must have
$\Diff_{b_1}\Diff_{b_2}\cdots \Diff_{b_N}f=0$. We will show
(Corollary~\ref{cor:abeltorsfree}) that this condition is not only
necessary but also sufficient.

\par
We note that this characterization easily implies that $f$ has a
periodic decomposition (with unprescribed periods) if and only if
$\Diff_{b_1}\Diff_{b_2}\cdots \Diff_{b_N}f=0$ for some pairwise
incommensurable $b_1,\dots,b_N$ and some $N\in\NN$.  This result was
already proved by Mortola and Peirone \cite{mortola/peirone:1999}.

\par
In this paper, we will consider a rather general situation: not only
translations on $\RR$ but mappings on arbitrary nonempty sets. So
the precise framework is the following.  We take an arbitrary
nonempty set $A$, and consider transformations $T:A\to A$. To such a
mapping we also associate a \emph{difference operator}
\begin{equation*}
\Diff_{T}f:=f\circ T-f.
\end{equation*}
A function is called then \emph{$T$-periodic} (or
$T$-\emph{invariant}), if $\Diff_{T}f=0$. This terminology is
naturally motivated by the case when $A=\RR$ and the transformation
$T$ is simply a translation by an element $a$ of $\RR$, i.e.,
$T(x):=T_a(x):=x+a$ for all $x\in\RR$. Note that in this case
$T$-periodicity of a function coincides with the usual notion of
$a$-periodicity.

\par Consider now  pairwise commuting transformations $T_1,T_2,\dots,
T_n:A\to A$. We say that a function $f:A\to \RR$ has a
$(T_1,T_2,\dots,T_n)$-\emph{periodic (or invariant) decomposition},
if
\begin{equation}\label{eq:adecomp}
f=f_1+\cdots+f_n\quad\mbox{with $f_j$ being $T_j$-periodic (i.e.,
$\Diff_{T_j}f_j=0$) for $j=1,\dots,n$}.
\end{equation}

We are looking for necessary and sufficient conditions which ensure
that such a periodic decomposition exists. A necessary condition is
once again clear as noted at the beginning: if
$f=f_1+f_2+\cdots+f_n$ holds with $f_j$ being $T_j$-periodic for
$j=1,\dots,n$, then using the commutativity of the transformations
we obtain
\begin{equation}\label{eq:diffn}
\Diff_{T_1}\Diff_{T_2}\cdots \Diff_{T_n}f=0+0+\cdots+0=0.
\end{equation}
As shown above, if we take $T_1=T_2$=translation by $a$ on $\RR$, we
see that having \eqref{eq:diffn} for some function $f:A\to\RR$, does
\emph{not} suffice for the existence of a periodic decomposition
\eqref{eq:adecomp}. (We remark that the already mentioned result
that $B(\R)$ has the decomposition property was in fact proved in
\cite{laczkovich/revesz:1990} by showing that the space $B(A)$ of
bounded functions on \emph{any} set $A$ has the decomposition
property with respect to \emph{arbitrary} commuting transformations
$T_1,T_2,\dots, T_n$.)

\par The decomposition problem for translation operators on $\RR$
originates  from  I.~Z.\ Ruzsa. He showed that the identity function
$\Id(x)=x$ can be decomposed into a sum of $a$- and $b$-periodic
functions, whenever $a/b$ is irrational. M.~Wierdl
\cite{wierdl:1984} extended this by showing that if
$a_1,\dots,a_n\in\RR$ are \emph{linearly independent} over $\QQ$ and
a function $f:\RR\to\RR$ satisfies \eqref{eq:diffeq}, then it has a
decomposition $f=f_1+f_2+\cdots+f_n$ with $a_j$-periodic functions
$f_j:\RR\to\RR$, $j=1,\dots,n$.

\par The periodic decomposition --  or invariant decomposition -- problem for \emph{arbitrary}, ``abstract''
transformations completely without structural restrictions on the
underlying set or on the function class was addressed in
\cite{farkas/revesz:2006}. There a certain Condition
(\textasteriskcentered) was presented, which was shown to be
necessary for the existence of periodic decompositions, moreover, it
was also proved to be sufficient for $n=1,2,3$ transformations. Now,
we restrict ourselves to the case of \emph{invertible}
transformations. Then the above mentioned Condition
(\textasteriskcentered) of \cite{farkas/revesz:2006} simplifies to
Condition \ref{cond:diffmod} of the present note, and we obtain both
necessity and sufficiency for any number $n\in\NN$ of
transformations.

We also reformulate the result in a particular case: we investigate
periodic, i.e., translation invariant decompositions of functions
defined on Abelian groups (Corollaries~\ref{cor:abel} and
\ref{cor:abeltorsfree}).

 In Section~\ref{appl}, we extend our result to
functions between torsion free Abelian groups (Theorem~\ref{tf2tf}).
As an immediate corollary of this extension, we obtain a positive
answer (Corollary~\ref{cor:integer}) to the question whether the
existence of a real valued periodic decomposition of an integer
valued function on $\RR$ (or more generally, on a torsion free
Abelian group) implies the existence of an integer valued periodic
decomposition with the same periods. This question was asked in
\cite{keleti/ruzsa:2007}, where several partial results have been
proved, for example that for functions defined on $\Z$ the analogous
 result holds. Here
we also show that this assertion about integer valued decompositions
is not true for functions on an arbitrary Abelian group.

A combination of our results with \cite{keleti:measurable}, yields a
characterization (Corollary~\ref{cor:measurable}) of those periods
$a_1,\ldots,a_n\in\R$ for which the existence of a real valued
measurable periodic decomposition of an $\R\to\Z$ function implies
the existence of an integer valued measurable periodic decomposition
with the same periods.

\section{Characterization of existence of periodic decompositions}

 Let $\mathcal{G}$ be the Abelian group generated by the
commuting, invertible transformations $T_1,\dots,T_n$ acting on a
set $A$, i.e., $\mathcal{G}:=\langle T_1,T_2,\dots, T_n\rangle$. For
an $x\in A$ we call the set $\{S(x):S\in\mathcal{G}\}$ the
\emph{orbit of $x$ under
 $\mathcal{G}$}. Such a set is often called an orbit of $\mathcal{G}$ as well, whereas this terminology shall not cause ambiguity.
 For a transformation $T:A\to A$ the orbits are understood as the orbits of $\langle T\rangle$.

The following observation helps to simplify the later arguments
considerably.
 Note that $T_j$-periodicity, hence the existence of a
$(T_1,T_2,\dots,T_n)$-periodic decomposition of a function, as well
as the validity of conditions involving difference operators (e.g.,
as in \eqref{eq:diffn}) is decided on the orbits of $\mathcal{G}$.
This means that we can always restrict considerations to the orbits
of $\mathcal{G}$.

\vskip1ex\par\noindent Now, using the following notation, we can
formulate the condition characterizing existence of periodic
decompositions and we can state our main result.

\begin{notation}\label{not:1}
For an Abelian group $\mathcal{H}$ and
$B=\{b_1,\dots,b_k\}\subseteq\mathcal{H}$ a nonempty set, we set
$[B]:=\bigcap_{i=1}^k \langle b_i\rangle$.
\end{notation}

\begin{condition}\label{cond:diffmod}For all orbits $\mathcal{O}$ of
$\mathcal{G}$, for all partitions
\begin{equation*}B_1\cup B_2\cup\cdots
\cup
B_N=\{T_1\rest_\mathcal{O},T_2\rest_\mathcal{O},\dots,T_n\rest_\mathcal{O}\}
\end{equation*}
and any element $S_j\in [B_j]$, $j=1,\dots, N$, we have that
\begin{equation}\label{eq:diffmod}
\Diff_{{S_1}}\dots
\Diff_{{S_N}}f\rest_\mathcal{O}=0\quad\mbox{holds}.
\end{equation}
\end{condition}

\begin{theorem}\label{th:sufficiency}
Let $T_1,\dots,T_n$ be pairwise commuting invertible transformations
on a set $A$. Let $f:A\to \RR$ be any function. Then $f$ has a
$(T_1,T_2,\dots, T_n)$-periodic decomposition if and only if it
satisfies Condition~\ref{cond:diffmod}.
\end{theorem}

\begin{remark}\hfill\rule{0pt}{0pt}\label{rem:egy}
\begin{enumerate}[1)]
\item Of course, $[B_j]$ is a cyclic group here,
and in the above Condition~\ref{cond:diffmod} it suffices to
consider only one of the \emph{generators} of $[B_j]$ instead of
\emph{all} elements $S_j\in [B_j]$. If $\mathcal{G}$ is torsion free
then there is a unique generator $S_j$ (up to taking possibly the
inverse), so Condition~\ref{cond:diffmod} simplifies.
\item If $T_k=T_{a_k}$ are translations on $\RR$, and if for some
$j$ there are incommensurable elements $a_k, a_m$ with the
corresponding transformations $T_{a_k}, T_{a_m}$ belonging to $B_j$,
then $[B_j]=\{\Id\}$. So for such partitions \eqref{eq:diffmod}
trivializes, hence it suffices to state Condition~\ref{cond:diffmod}
for partitions $B_1\cup B_2\cup\cdots \cup B_N$ for which the
elements  belonging to each $B_j$, $j=1,\dots,N$ are all
commensurable. According to the above it also suffices to consider
$S_j$ to be the translation by the least common multiple $b_j$ of
the elements in $B_j$.
\end{enumerate}
\end{remark}

Assume now that the underlying set is a group $A=\mathcal{G}$ and
the transformations $T_j$ are one-sided, say right, multiplications
by elements $a_j\in \mathcal{G}$. Then commutativity of the
transformations $T_j$ is equivalent to commutativity of the
generating elements $a_j$. Note that in this case the orbits of the
group are the left cosets of the subgroup generated by
$\{a_1,\dots,a_n\}$ in $\mathcal{G}$. Then each transformation acts
on each orbit in the same way, so in Condition~\ref{cond:diffmod} we
do not have to restrict  anything to the orbits. Hence
Theorem~\ref{th:sufficiency} gives the following.

\begin{corollary}
\label{cor:group} Let $\mathcal{G}$ be a group and $a_1,\dots,a_n\in
\mathcal{G}$ commuting pairwise with each other. Then a function
$f:\mathcal{G}\to \RR$ decomposes into a sum of
right-$a_j$-invariant functions, $f=f_1+\cdots+ f_n$, if and only if
for all partitions $B_1\cup B_2\cup\cdots \cup
B_N=\{a_1,\dots,a_n\}$ and for any element $b_j\in [B_j]$ (see
Notation \ref{not:1}) one has
\begin{equation*}
\Diff^{(r)}_{{b_1}}\dots \Diff^{(r)}_{{b_N}}f=0,
\end{equation*}
with $\Diff^{(r)}_{{a}}$ denoting the right difference operator:
$\bigl(\Diff^{(r)}_{{a}}\bigr)f(x)=f(xa)-f(x)$.
\end{corollary}

In the special case, when $T_j$ are translations on an Abelian group
written additively, the above yields immediately:

\begin{corollary} \label{cor:abel}
Let $\mathcal{A}$ be an (additive) Abelian group and
$a_1,\dots,a_n\in \mathcal{A}$. A function $f:\mathcal{A}\to \RR$
decomposes into a sum of $a_j$-periodic functions, $f=f_1+\cdots+
f_n$, if and only if for all partitions $B_1\cup B_2\cup\cdots \cup
B_N=\{a_1,\dots,a_n\}$ and for any element $b_j\in [B_j]$ (see
Notation \ref{not:1}) one has
\begin{equation*}
\Diff_{{b_1}}\dots \Diff_{{b_N}}f=0.
\end{equation*}
\end{corollary}

\noindent In a torsion free Abelian group $\mathcal{A}$, we call the
unique generator $b$ (up to taking possibly the inverse) of the
cyclic group $\langle b_1\rangle\cap\dots\cap\langle
b_m\rangle=[\{b_1,b_2,...,b_m\}]$ the \emph{least common multiple}
of the elements $b_1,b_2,\dots, b_m\in\mathcal{A}$. Note that with
this terminology we have for example that the least common multiple
of $1$ and $\sqrt{2}$ in the group $(\R,+)$ is $0$. This is also the
case in general: $a$ and $b$ are by definition incommensurable in
$\mathcal{A}$, if  $[\{a,b\}]=\{e\}$, i.e., if their least common
multiple is the unit element $e\in\mathcal{A}$.

We can now reformulate Theorem \ref{th:sufficiency} in this special
case as follows (see also Remark \ref{rem:egy}).
\begin{corollary} \label{cor:abeltorsfree}
Let $\mathcal{A}$ be a torsion free Abelian group and
$a_1,\dots,a_n\in \mathcal{A}$. A function $f:\mathcal{A}\to \RR$
decomposes into a sum of $a_j$-periodic functions, $f=f_1+\cdots+
f_n$, if and only if for all partitions $B_1\cup B_2\cup\cdots \cup
B_N=\{a_1,\dots,a_n\}$ and $b_j$ being the least common multiple of
the elements in $B_j$ one has
\begin{equation*}
\Diff_{{b_1}}\dots \Diff_{{b_N}}f=0.
\end{equation*}
\end{corollary}

Theorem~\ref{th:sufficiency} can be also applied if $\mathcal{G}$ is
a non-Abelian group and among the transformations there are both
\emph{left} multiplications by certain pairwise commuting elements,
and some \emph{right} multiplications by certain further elements
again pairwise commuting among themselves. Indeed, a left and a
right multiplication can always be interchanged in view of the
associativity law, so this way we get pairwise commuting invertible
transformations, and Theorem~\ref{th:sufficiency} can be applied.
However, in this case the orbits are double cosets on which our
transformations can act differently, so in this case one cannot
simplify Condition~\ref{cond:diffmod} as in
Corollary~\ref{cor:group}.

 We also remark that if we look for continuous
decomposition of continuous functions $f:\RR\to\RR$, then
Condition~\ref{cond:diffmod} is again insufficient: if $f(x)=x$ and
$a_1/a_2\not\in \QQ$ then Condition~\ref{cond:diffmod} is satisfied
but $f(x)=x$ is not a sum of two \emph{continuous}, periodic
functions, because it is not bounded. In fact, $f(x)=x$ is not even
the sum of two \emph{measurable} periodic functions (see, e.g., in
\cite{laczkovich/revesz:1990}), so Condition~\ref{cond:diffmod} is
insufficient even for measurable decomposition of continuous
functions.

\section{The proof of Theorem~\ref{th:sufficiency}}

\par The necessity of Condition~\ref{cond:diffmod} can already be found in
\cite[Theorem 4]{farkas/revesz:2006} and is proved there even for
not necessarily invertible transformations. We give a
 straightforward proof for the reader's
convenience.
\begin{proposition}[{\textbf{Necessity}}]\label{th:necessity}
Suppose that $T_1,\dots,T_n$ are pairwise commuting invertible
 transformations on a set $A$. If $f:A\to
\RR$ has a $(T_1,T_2,\dots, T_n)$-periodic decomposition then
Condition~\ref{cond:diffmod} is satisfied.
\end{proposition}
\begin{remark}
Actually, the proof below yields mutatis mutandis the analogous
result for functions $f:A\to \Gamma$, where $\Gamma$ is an arbitrary
Abelian group, written additively.
\end{remark}
\begin{proof}
We can assume that the group of transformations $\langle T_1,\dots,
T_n\rangle$ (generated by $T_1,\dots, T_n$) acts on $A$
transitively, i.e., that $A$ is already one orbit under the action
of the transformations.  Given a partition $B_1\cup B_2\cup\cdots
\cup B_N= \{T_1,T_2,\dots,T_n\}$ and $S_j\in [B_j]$, we have to show
$\Diff_{{S_1}}\dots \Diff_{{S_N}}f=0$.

 Note that for $T_i\in B_j$
the function $f_i$ is $S_j$-periodic as $S_j=T_i^{m_i}$ for some
$m_i$ (without loss of generality we can assume $m_i\geq 0$,
otherwise we could repeat the argument for $S_j^{-1}=T_i^{-m_i}$),
and

\begin{equation*}
f_i(x)=f_i\bigl(T_i(x)\bigr)=f_i\bigl(T_i^2(x)\bigr)=\cdots=
f_i\bigl(T_i^{m_i}(x)\bigr)=f_i\bigl(S_j(x)\bigr).
\end{equation*}
So by summing up for each fixed $j$ the functions $f_i$
corresponding to those transformations $T_i$ which belong to $B_j$,
we get the functions $g_j:=\sum_{T_i\in B_j} f_i$ which are still
$S_j$-periodic, therefore $f=g_1+\cdots+g_N$ is an $(S_1,S_2,\dots,
S_N)$-periodic decomposition of $f$. Hence, as $S_1,S_2,\dots,S_N$
are pairwise commuting and so are $\Diff_{S_1},\Diff_{S_2},\dots,
\Diff_{S_N}$, we obtain that indeed $\Diff_{{S_1}}\Diff_{{S_2}}\dots
\Diff_{{S_N}}f=0$.
\end{proof}

To prove sufficiency of Condition~\ref{cond:diffmod} we will need
the following lemma, which is a slightly modified version of a
result from \cite{farkas/revesz:2006}, where the same result was
proved for real valued functions and the invertibility of the
transformations were not assumed. For the sake of completeness we
present the same proof in a compacter form for our case.

\begin{lemma}\label{lem:diffinvsolv} Let $T,S$ be commuting invertible transformations of $A$ and let
$G:A\to\Gamma$ be a function with values in the (additive) Abelian
group $\Gamma$ and satisfying $\Diff_T G=0$. Then there exists a
function $g: A\to \Gamma$ satisfying both $\Diff_T g= 0$ and
$\Diff_S g =G$ if and only if
\begin{equation}\label{eq:twotransfcondi}
\sum_{i=0}^{n-1} G(S^i({x})) =0\quad\mbox{whenever $T^m
({x})=S^n({x})$ for some  $m\in\ZZ$, $n\in \NN$, ${x}\in A$}.
\end{equation}
\end{lemma}
\begin{proof}
Suppose first that there exists a $T$-periodic $g:A\to \Gamma$ with
$\Diff_S g=G$, and also that $T^m({x})=S^n({x})$ for some $m\in
\ZZ$, $n\in\NN$ and ${x}\in A$. Then
\begin{equation*}
\sum_{i=0}^{n-1} G(S^i({x}))
=\sum_{i=0}^{n-1}\bigl(g(S^iS({x}))-g(S^i({x})\bigr)
=g(S^n({x}))-g({x})=g(T^m({x}))-g({x})=0,
\end{equation*}
by the $T$-periodicity of $g$. So condition
\eqref{eq:twotransfcondi} is necessary.

\par We now prove the sufficiency of this condition. Let us consider the set
$\widetilde{A}$ of all orbits of the cyclic group $\langle T
\rangle$.
 Since $G$ is $T$-periodic, it is constant on each orbit $\xelm{x}\in
\widetilde{A}$, i.e., $G({x})=G({x}')$ if ${x},{x}'\in \xelm{x}$
(with a small abuse of notation we will write $\xelm{x}$ for the
orbit of ${x}$). So the function $\widetilde{G}:\widetilde{A}\to
\Gamma$ defined by $\widetilde{G}(\xelm{x})=G(x)$ is well-defined.
Because of commutativity the transformation $S$ maps orbits of
$\langle T\rangle$ into orbits, hence we can define
$\widetilde{S}:\widetilde{A}\to \widetilde{A}$ by
$\widetilde{S}(\xelm{x}):=\xelm{y}$ with $\xelm{y}$ the orbit of
$S({x})$. Now we pass to $\widetilde{A}$, and notice that
\eqref{eq:twotransfcondi} implies
\begin{equation}\label{eq:onetransfcondi}
\sum_{i=0}^{n-1} \widetilde{G}(\widetilde{S}^i(\xelm{x}))
=0\quad\mbox{whenever $\widetilde{S}^n (\xelm{x})=\xelm{x}$ for some
$n\in \NN$ and $\xelm{x}\in \widetilde{A}$}.
\end{equation}

Consider the orbits of $\langle \widetilde{S}\rangle$ in
$\widetilde{A}$. By the axiom of choice we select for each such
orbit ${\beta}\subset \widetilde{A}$ an element $\xelm{x}_{\beta}\in
{\beta}$. We claim that the function defined as follows
(understanding empty sums as $0$) is well defined:
\begin{equation*}
\widetilde{g}(\xelm{x}):=
\begin{cases}\displaystyle
\widetilde{G}(\xelm{x}_{{\beta}})-\sum_{i=0}^{n-1}
\widetilde{G}(\widetilde{S}^i(\xelm{x})),&\mbox{if $\xelm{x}\in
{\beta}$ and $\widetilde{S}^n(\xelm{x})=\xelm{x}_{{\beta}}$
with $n\geq 0$},\\[2ex]\displaystyle
\widetilde{G}(\xelm{x}_{{\beta}})+\sum_{i=n}^{-1}
\widetilde{G}(\widetilde{S}^i(\xelm{x})),&\mbox{if $\xelm{x}\in
{\beta}$ and $\widetilde{S}^n(\xelm{x})=\xelm{x}_{{\beta}}$ with
$n<0$}.
\end{cases}
\end{equation*}
Indeed, if both $\widetilde{S}^n(\xelm{x})=\xelm{x}_{{\beta}}$ and
$\widetilde{S}^m(\xelm{x})=\xelm{x}_{{\beta}}$ hold with, say,
$n-m>0$, then
$\widetilde{S}^{n-m}\widetilde{S}^m(\xelm{x})=\widetilde{S}^m(\xelm{x})$,
so \eqref{eq:onetransfcondi} yields that the difference of the two
different expressions that define $\widetilde{g}(\xelm{x})$ is
\begin{equation*}
\sum_{i=m}^{n-1}
\widetilde{G}(\widetilde{S}^i(\xelm{x}))=\sum_{i=0}^{n-m-1}
\widetilde{G}(\widetilde{S}^i\widetilde{S}^m(\xelm{x}))=0.
\end{equation*}
For $\widetilde{S}^n(\xelm{x})=\xelm{x}_{{\beta}}$, $n> 0$, we have
$\widetilde{S}^{n-1}\widetilde{S}(\xelm{x})=\xelm{x}_{{\beta}}$ and
hence, by definition of $\widetilde{g}$,
\begin{align*}
\Diff_{\widetilde{S}}\widetilde{g}(\xelm{x})&=\widetilde{g}(\widetilde{S}(\xelm{x}))-\widetilde{g}(\xelm{x})\\
&=\left(\widetilde{G}(\xelm{x}_{{\beta}})
-\sum_{i=0}^{n-2}\widetilde{G}(\widetilde{S}^i\widetilde{S}(\xelm{x}))\right)
-\left(\widetilde{G}(\xelm{x}_{{\beta}})
-\sum_{i=0}^{n-1}\widetilde{G}(\widetilde{S}^i(\xelm{x}))\right)
=\widetilde{G}(\xelm{x}).
\end{align*}
It follows similarly
$\Diff_{\widetilde{S}}\widetilde{g}(\xelm{x})=\widetilde{G}(\xelm{x})$
also in the cases $n=0$, $n<0$.

\par Now we pull back
$\widetilde{g}:\widetilde{A}\to\widetilde{A}$ to $A$, that is we set
$g({x}):=\widetilde{g}(\xelm{x})$, where $\xelm{x}$ is the orbit of
${x}$ under $\langle T\rangle$. It is straightforward that $\Diff_T
g=0$ and $\Diff_S g=G$.
\end{proof}

\noindent We complete the proof of Theorem~\ref{th:sufficiency} by
proving the sufficiency of Condition~\ref{cond:diffmod}.

\begin{remark}\label{rem:N} The reader will have no difficulty to
verify that the proof below yields also the following assertions.
\begin{enumerate}[1)]
\item
Theorem \ref{th:sufficiency} remains valid  for functions taking
values in an arbitrary divisible torsion free group (which can be
written in the form $\bigoplus_{I} \QQ = \QQ^{I}$ for some index set
$I$).
\item There is a positive integer $M=M(T_1,\dots,T_n)$
such that, whenever $f$ takes its values in an additive subgroup
$\Gamma$ of $\QQ^{I}$, then the functions in the periodic
decomposition can be chosen to have values $\frac xM$ with $x\in
\Gamma$.
\item  If $T_1$ is of infinite order, $T_1^{m_1}=T_2^{m_2}=\ldots=T_k^{m_k}$
($m_1\in\N$, $m_2\ldots,m_k\in\Z\sm\{0\}$) and $T_{k+1},\ldots,T_n$
are unrelated to $T_1$ (see the definition below), then as the
constant $M(T_1,T_2,\dots, T_n)$ we can take $m_1 M(T_2,\dots,T_n)
M(T_{k+1},\dots,T_n)$.
\end{enumerate}
\end{remark}

\begin{proof}[Proof of the sufficiency of Condition~\ref{cond:diffmod} in
Theorem~\ref{th:sufficiency}] Our proof is by induction on the
number $n\in\NN$ of transformations. The case $n=1$ is obvious.

As pointed out before, we assume without loss of generality that
$\mathcal{G}:=\langle T_1,\dots, T_n\rangle$ acts transitively on
$A$, i.e., $A$ is one orbit of $\mathcal{G}$.  Following
\cite{farkas/revesz:2006} we say that the transformations $T_i$ and
$T_j$ are \emph{related}, if there are $m,k\in \ZZ\setminus\{0\}$
with $T_i^m=T_j^k$. This is clearly an equivalence relation.

If  possible, let us take as $T_1$ an element of infinite order from
$\{T_1,T_2,\dots, T_n\}$.  Let us also assume for notational
convenience that $\{T_1,T_2,\dots T_k\}$, $k\leq n$, are exactly the
elements among the $T_j$s being related to $T_1$. This means
particularly that there exist $m_j\in\ZZ\setminus\{0\}$, $m_1\in
\NN$ such that $T_j^{m_j}= T_1^{m_1}=:U_1$ for all $j=1,\dots, k$,
and the other elements $\{T_{k+1},\dots, T_n\}$ are then all
unrelated to $T_j$, $j=1,\dots, k$. Note that if all the elements
$T_1, T_2,\dots, T_n$ have finite order, then they are also all
related, so $n=k$.

 We define the functions $g:=\Diff_{T_1} f$,
$h:=\Diff_{U_1} f$. In case $U_1=T_1$, these two functions coincide,
but the following arguments still remain valid. It might also happen
that $U_1=\Id$, in this case $h=0$, but neither does this effect the
validity of the following.

We will apply the induction hypothesis to $g$ and $h$. For this
purpose we first check that Condition~\ref{cond:diffmod} is
satisfied for the function $g$ and the transformations
$\{T_2,T_3,\dots,T_n\}$ and for the function $h$ and the
transformations $\{T_{k+1},\dots, T_n\}$, respectively. By
assumption, we have Condition~\ref{cond:diffmod} for $f$ and the
transformations $\{T_1,T_2,\dots, T_n\}$; we are to apply this by
choosing the partitions $B_1\cup B_2\cup \dots\cup B_N$ in a
particular way.

Considering partitions of $\{T_1,T_2,\dots,T_n\}$ with
$B_1:=\{T_1\}$ and the other blocks being arbitrary and taking
$S_1=T_1\in[B_1]$ and $S_j\in [B_j]$ arbitrary, we see that
Condition~\ref{cond:diffmod} is satisfied for $g$ and
transformations $\{T_2,T_3,\dots,T_n\}$. Similarly, if we consider
$B_1=\{T_1,T_2,\dots,T_k\}$, $S_1=U_1$ and the other blocks
arbitrary, we see that $h$ satisfies Condition~\ref{cond:diffmod}
with transformations $\{T_{k+1},T_{k+2},\dots, T_n\}$.

Now the inductive hypothesis yields the two decompositions
\begin{align*}
g&=g_2+\cdots+g_k+g_{k+1}+\cdots+g_n,&\qquad \mbox{with }\Diff_{T_j}
g_j=0, &\quad 2\leq j\leq n,\\
\mbox{and}\quad  h&=h_{k+1}+\cdots+h_n,&\qquad \mbox{with
}\Diff_{T_j} h_j=0, &\quad k+1\leq j\leq n.
\end{align*}
(Note again that, if incidentally $k=n$, then $h=0$ by assumption.)
For all $j=2,\dots, n$ we define the function
\begin{equation*}
G_j({x}):=\frac{1}{m_1}\sum_{\mu=0}^{m_1-1} g_j(T_1^\mu ({x})),
\end{equation*}
 which
is, of course, $T_j$-periodic. Moreover, for $j=2,\dots,k$ we can
even claim $\Diff_{T_1} G_j=0$. Indeed, one has
\begin{equation*}
\Diff_{T_1}G_j({x})=G_j(T_1({x}))-G_j({x})=\frac{1}{m_1}\bigl(g_j(T_1^{m_1}
({x}))-g_j({x})\bigr)=0,
\end{equation*}
because $T_1^{m_1}=T_j^{m_j}$ and $g_j$ is $T_j$-periodic. By
definition we can write,
\begin{equation*}
h({x})=f(T_1^{m_1}(x))-f({x})=\sum_{\mu=0}^{m_1-1} g(T_1^\mu({x})) =
\sum_{\mu=0}^{m_1-1} \sum_{j=2}^n g_j(T_1^\mu({x}))= m_1\sum_{j=2}^n
G_j({x}),
\end{equation*}
whence the decomposition of $h$ entails
\begin{equation}\label{eq:twodecomp}
\frac{1}{m_1}\sum_{j=k+1}^n h_j({x}) = \sum_{j=2}^n G_j({x}).
\end{equation}
Now, the functions $F_j$, $j=2,\dots,n$, defined by
\begin{equation*}
F_j:= \left\{%
\begin{array}{ll}
g_j-G_j, & 2\leq j\leq k, \\
\displaystyle g_j-G_j+\frac{1}{m_1}h_j, & k+1\leq j\leq n
\end{array}%
\right.
\end{equation*}
are undoubtedly $T_j$-periodic. According to \eqref{eq:twodecomp} we
still have
\begin{equation}\label{eq:Freprez}
\sum_{j=2}^n F_j = \sum_{j=2}^n g_j = g = \Diff_{T_1}f.
\end{equation}
Now we prove that we can apply Lemma \ref{lem:diffinvsolv} with
$S=T_1$ and $T=T_j$ to all functions $F_j$.  For the indices $j\geq
k+1$ these transformations are unrelated to $T_1$. This means that
$T_1^m=T_j^{m'}$ can not hold for $m,m'\in \ZZ\setminus\{0\}$. Nor
is it possible that $T_1^m=\Id=T_j^0$ with $m\in\ZZ\setminus\{0\}$,
because $T_1$ was chosen to be of infinite order. (Should this
choice be impossible, then all elements are related, i.e., $n=k$ and
this case is empty.) So we see that for $j\geq k+1$ condition
\eqref{eq:twotransfcondi} of Lemma \ref{lem:diffinvsolv} is void,
whence the existence of a ``lift-up'' $f_j$ with $\Diff_{T_1}
f_j=F_j$ and $\Diff_{T_j} f_j=0$ is immediate. Let us consider the
cases of $j=2,\dots,k$. Then the two transformations $T=T_1$ and
$S=T_j$ are related. So let now $T_j^{n_j}=T_1^m$ for some $m\in\NN$
and $n_j\in\ZZ\setminus\{0\}$. Let us take now $k_j:=\min \{\ell\in
\NN~:~ \exists \nu\in\ZZ\setminus\{0\}\;T_j^\nu=T_1^\ell\}$,
$j=2,\dots,k$. Clearly, we have then $k_j|m$ and $k_j|m_1$. From
this and the $T_j$-periodicity of $g_j$ we obtain
\begin{equation*}
\frac{1}{m}\sum_{\mu=0}^{m-1}
g_j(T_1^\mu({x}))=\frac{1}{k_j}\sum_{\mu=0}^{k_j-1} g_j(T_1^\mu(
{x}))=\frac{1}{m_1}\sum_{\mu=0}^{m_1-1} g_j(T_1^\mu( {x}))=G_j({x}).
\end{equation*}
Therefore, using also $\Diff_{T_1} G_j({x})=0$ for $j=2,\dots,k$, we
get
\begin{equation*}
\sum_{\mu=0}^{m-1} F_j(T_1^\mu( {x})) = \sum_{\mu=0}^{m-1}
g_j(T_1^\mu({x})) - \sum_{\mu=0}^{m-1} G_j(T_1^\mu({x})) = m \cdot
G_j({x})- \sum_{\mu=0}^{m-1} G_j({x}) = 0.
\end{equation*}
This shows that for $T=T_j$ and $S=T_1$ the assumptions of Lemma
\ref{lem:diffinvsolv} are satisfied, hence the application of this
lemma furnishes $T_j$-periodic functions $f_j:A\to \RR$ with
$\Diff_{T_1}f_j=F_j$, $j=2,\dots,k$.

Finally, we set $f_1:=f-(f_2+f_3+\cdots +f_n)$. Using
\eqref{eq:Freprez} we see that
\begin{equation*}
\Diff_{T_1}f_1=\Diff_{T_1}f-(\Diff_{T_1}f_2+\Diff_{T_1}f_3+\cdots+\Diff_{T_1}f_n)=\Diff_{T_1}f-(F_2+F_3+\cdots+F_n)=0.
\end{equation*}
Thus $f=f_1+f_2+\cdots+f_n$ is a desired periodic decomposition of
$f$.
\end{proof}

\section{Generalizations and applications}\label{appl}

We start this section by proving that Corollary
\ref{cor:abeltorsfree} also holds for $\Gamma$-valued functions for
arbitrary torsion free Abelian group $\Gamma$.
\begin{theorem}\label{tf2tf}
Let $\mathcal{A}, \Gamma$ be torsion free Abelian groups and
$a_1,\dots,a_n\in \mathcal{A}$. A function $f:\mathcal{A}\to \Gamma$
decomposes into a sum of $a_j$-periodic functions
$f_j:\mathcal{A}\to \Gamma$, $f=f_1+\cdots+f_n$, if and only if for
all partitions $B_1\cup B_2\cup\cdots \cup B_N=\{a_1,\dots,a_n\}$
and $b_j$ being the least common multiple of the elements in $B_j$
one has
\begin{equation*}
\Diff_{{b_1}}\dots \Diff_{{b_N}}f=0.
\end{equation*}
\end{theorem}
\begin{proof} Since the proof of the necessity worked for
arbitrary $\IA, \Gamma$, it suffices to prove sufficiency. This we
prove by induction on $n$, the case $n=1$ being obvious. So let
$n>1$ and let us assume that the assertion holds for any smaller
number of translates.

We can (re)number the translates $\{a_2,\ldots,a_n\}$ so that
$a_2,\ldots,a_k$ are exactly the elements that have common multiple
with $a_1$. Let $b$ denote the least common multiple of
$a_1,\ldots,a_k$ and take the integers $m_1,\dots,m_k$ such that
$b=m_1 a_1=\cdots=m_k a_k$. Using that $\mathcal{A}$ is torsion
free, it is easily seen that $m_1,\dots,m_k$ are relatively primes
and so there exist integers $d_1,\dots,d_k$ for which $d_1 m_1 +
\cdots + d_k m_k = 1$. We can suppose that $\Gamma \leq \QQ^{I}$,
since a torsion free Abelian group can be embedded as a subgroup in
$\QQ^{I}$ for some $I$. As we saw in Remark \ref{rem:N}.2,  for a
function $f$ with values in $\Gamma \leq \QQ^{I}$ the proof of
Theorem \ref{th:sufficiency} gives only a decomposition with values
in $\frac{1}{M} \Gamma$ where $M = M(T_1,T_2,\dots,T_n)$ is a
certain positive integer. Furthermore, by Remark \ref{rem:N}.3,
 the constant $m_1
M(T_2,\dots,T_n) M(T_{k+1},\dots,T_n)$ is an appropriate choice for
$M(T_1,T_2,\dots,T_n)$. By the induction hypothesis, however, we can
already assume that $M(T_2,\dots,T_n) = M(T_{k+1},\dots,T_n) = 1$,
so we get a decomposition $f = f^1_1 + \cdots + f^1_n$ where $f^1_j:
\IA \to \frac{1}{m_1}\Gamma$ is an $a_j$-periodic function
($j=1,\dots,n$). Repeating the same argument for $a_i$ ($1 \leq i
\leq k$) instead of $a_1$, we obtain $k$ decompositions of $f$ with
the following properties:
\begin{equation*}
f = f^i_1 + f^i_2 + \cdots + f^i_n \mbox{, where $f^i_j : \IA \to
\frac{1}{m_i}\Gamma$ is $a_j$-periodic } (i=1,\dots, k).
\end{equation*}
\noindent Take the linear combination of the decompositions $f^i_j$
with coefficients $d_i m_i$:
\begin{equation*}
f_j := \sum_{i=1}^{k} d_i m_i f^i_j.
\end{equation*}
Being the linear combination of $a_j$-periodic functions, $f_j$ is
$a_j$-periodic, too. Moreover, $f_j$ is $\Gamma$-valued because the
functions $m_i f^i_j$ are all $\Gamma$-valued ($i=1,\dots, k$).
Finally, the functions $f_j$ indeed give a suitable decomposition of
$f$, because
\begin{equation*}
\sum_{j=1}^{n} f_j = \sum_{j=1}^{n} \sum_{i=1}^{k} d_i m_i f^i_j =
\sum_{i=1}^{k} d_i m_i \sum_{j=1}^{n} f^i_j = \left( \sum_{i=1}^{k}
d_i m_i \right) \cdot f = f.
\end{equation*}
\end{proof}
\noindent  Let $\IA$ be a torsion free Abelian group. Using the
above theorem for $\Gamma = \RR$ and then for $\Gamma = \ZZ$, we get
that for a function $f : \IA \to \ZZ$ the existence of a real valued
periodic decomposition and the existence of an integer valued
periodic decomposition are both equivalent with the same condition.
Thus  we obtain the following.
\begin{corollary}\label{cor:integer}
If an integer valued function $f$ on a torsion free Abelian group
$\IA$ decomposes into the sum of $a_j$-periodic real valued
functions, $f=f_1+f_2+\cdots+f_n$ for some $a_1,\dots,a_n$, then $f$
also decomposes into the sum of $a_j$-periodic integer valued
functions, $f=g_1+g_2+\cdots+g_n$.
\end{corollary}

\noindent The following example shows that neither
Theorem~\ref{tf2tf} nor Corollary~\ref{cor:integer} holds for every
Abelian group $\IA$. For instance let $\IA$ be $\ZZ_2 \times \ZZ_2$
and $a_1 = (1,0)$; $a_2 = (0,1)$; $a_3 = (1,1) \in \ZZ_2 \times
\ZZ_2$. We can represent functions $f:\ZZ_2 \times \ZZ_2 \to \RR$
 by $2 \times 2$ real matrices: $\bigl(
\begin{smallmatrix}
f(0,0) & f(0,1) \\
f(1,0) & f(1,1)
\end{smallmatrix}
\bigr)$. Let us consider the following decomposition
\begin{equation*}
\begin{pmatrix}
0 & 1 \\
1 & 1
\end{pmatrix}=
\begin{pmatrix}
0 & \frac{1}{2}\\
0 & \frac{1}{2}
\end{pmatrix}
+
\begin{pmatrix}
0 & 0 \\
\frac{1}{2} & \frac{1}{2}
\end{pmatrix}
+
\begin{pmatrix}
0 & \frac{1}{2} \\
\frac{1}{2} & 0
\end{pmatrix}.
\end{equation*}
The functions $f_j : \ZZ_2 \times \ZZ_2 \to \RR$ represented by the
matrices on the right are clearly $a_j$-periodic respectively for
$j=1,2,3$. However, there does not exist an integer valued
decomposition since the sum of the entries of a matrix that
represents an $a_j$-periodic integer valued function $g_j : \ZZ_2
\times \ZZ_2 \to \ZZ$ must be even ($j=1,2,3$). Also the sums of
such matrices do have this property, whereas  $\bigl(
\begin{smallmatrix}
0 & 1 \\
1 & 1
\end{smallmatrix}
\bigr)$ does not.

\par Finally, for measurable decompositions we can
draw the following consequence of our results. Proposition 3.3 in
\cite{keleti:measurable} tells us that the equivalence of i) and ii)
in Corollary \ref{cor:measurable} below is valid, whenever the
statement of Corollary \ref{cor:integer} is true, what we  already
have proved. This consideration yields the following.
\begin{corollary}\label{cor:measurable} The following
two assertions are equivalent for arbitrary real numbers
$a_1,\dots,a_n$.
\begin{enumerate}[i)]
\item If an integer valued function $f:\RR \to \ZZ$ can be decomposed as
$f = f_1 + \cdots + f_n$ such that each $f_j$ is an $a_j$-periodic
measurable $\RR \to \RR$ function then $f$ can be also decomposed as
$f = g_1 + \cdots + g_n$ such that each $g_j$ is an $a_j$-periodic
integer valued measurable function.
\item If $B_1,\dots,B_N$ are the equivalence classes of $\{a_1,\dots,a_n\}$
with respect to the relation $a \sim b \Leftrightarrow a/b \in \QQ$,
and $b_j$ denotes the least common multiple of the numbers in $B_j$
($j=1,\dots,N$), then
$\frac{1}{b_1},\dots,\hskip-0.5pt\frac{1}{b_N}$\hskip-1pt are
linearly independent over $\QQ$.
\end{enumerate}
\end{corollary}

\parindent0pt


\begin{thebibliography}{99}

\bibitem{farkas/revesz:2006}
B. Farkas, Sz.Gy. R\'{e}v\'{e}sz,
\newblock \textit{Decomposition as the sum of invariant functions with respect to commuting transformations},
\newblock Aequationes Math., to appear.


\bibitem{gajda:1992}
Z. Gajda,
\newblock \textit{Note on decomposition of bounded functions into the sum of periodic terms},
\newblock Acta Math. Hung.  \textbf{59} (1992), no.~1-2, 103--106.


\bibitem{kadets/shumyatskiy:2000}
V.M. Kadets, S.B. Shumyatskiy,
\newblock \textit{Averaging Technique in the Periodic Decomposition Problem},
\newblock Mat. Fiz. Anal. Geom. \textbf{7} (2000), no.~2, 184--195.


\bibitem{kadets/shumyatskiy:2001}
V.M. Kadets, S.B. Shumyatskiy,
\newblock \textit{Additions to the Periodic Decomposition Theorem},
\newblock Acta Math. Hungar. \textbf{90} (2001), no.~4, 293--305.


\bibitem{keleti:1996}
T. Keleti,
\newblock\textit{Difference functions of periodic measurable functions},
\newblock PhD dissertation, ELTE, Budapest, 1996.


\bibitem{keleti:1997}
T. Keleti,
\newblock \textit{On the differences and sums of periodic measurable functions},
\newblock Acta Math. Hungar. \textbf{75}(1997), no.~4, 279--286.


\bibitem{keleti:1998}
T. Keleti,
\newblock\textit{Difference functions of periodic measurable functions},
\newblock Fund. Math.  \textbf{157} (1998), 15--32.

\bibitem{keleti:measurable}
T. Keleti, \textit{Periodic decomposition of measurable integer
valued functions}, J. Math. Anal. Appl., to appear.

\bibitem{keleti/ruzsa:2007}
Gy. K\'arolyi, T. Keleti, G. K\'os, I.Z. Ruzsa,
\newblock \textit{Periodic decomposition of integer valued functions},
\newblock Acta Math. Hungar., to appear.

\bibitem{laczkovich/revesz:1989}
M. Laczkovich, Sz.Gy. R\'{e}v\'{e}sz,
\newblock \textit{Periodic decompositions of continuous functions},
\newblock Acta Math. Hungar. \textbf{54}(1989), no.~3-4, 329--341.


\bibitem{laczkovich/revesz:1990}
M. Laczkovich, Sz.Gy. R\'{e}v\'{e}sz,
\newblock \textit{Decompositions into periodic functions belonging to a given Banach space},
\newblock Acta Math. Hung. \textbf{55}(3-4) (1990), 353--363.


\bibitem{mortola/peirone:1999}
S. Mortola, R. Peirone,
\newblock \textit{The sum of periodic functions}
\newblock Boll. Un. Mat. Ital. \textbf{8} 2-B (1999), 393--396.


\bibitem{wierdl:1984}
M. Wierdl,
\newblock \textit{Continuous functions that can be represented as the sum of finitely many periodic functions},
\newblock Mat. Lapok \textbf{32} (1984) 107--113 (in Hungarian).

\end{thebibliography}
\end{document}